\documentclass[a4paper,10pt,leqno]{article}

\usepackage[english]{babel}

\usepackage{amsfonts,amsmath}
\usepackage{amssymb}
\usepackage{eucal}
\usepackage{mathrsfs}
\usepackage[hypertexnames=false]{hyperref}
\usepackage{bbm}
\usepackage{graphicx}
\usepackage{caption}
\usepackage{subcaption}
\usepackage{authblk}

\usepackage[usenames,dvipsnames]{color}

\newcommand{\R}{\mathbb{R}}

\providecommand{\norm}[1]{\left\lVert#1\right\rVert}

\begin{document}

\thispagestyle{empty}

\title{A primal-dual approach for a total variation Wasserstein flow}
\author[1]{Martin Benning}
\author[2]{ Luca Calatroni}
\author[3]{Bertram D\"{u}ring}
\author[4]{Carola-Bibiane Sch\"{o}nlieb}
\affil[1]{\small{Magnetic Resonance Research Centre, University of Cambridge, UK}}
\affil[2]{Cambridge Centre for Analysis, University of Cambridge, UK}
\affil[3]{Department of Mathematics, University of Sussex, UK}
\affil[4]{Department of Applied Mathematics and Theoretical Physics, University of Cambridge, UK}

\date{}
\maketitle


\begin{abstract}
We consider a nonlinear fourth-order diffusion equation that arises in denoising of image densities. We propose an implicit time-stepping scheme that employs a primal-dual method for computing the subgradient of the total variation seminorm. The constraint on the dual variable is relaxed by adding a \emph{penalty term}, depending on a parameter that determines the weight of the penalisation. The paper is furnished with some numerical examples showing the denoising properties of the model considered.
\end{abstract}


\section{Introduction} \label{sec:int}

For an open and bounded domain $\Omega\subset\R^2$ with Lipschitz boundary, we consider the following equation for $u=u(t,x)$
\begin{equation} \label{TVWassflow}
\begin{aligned}
& u_t=\nabla\cdot(u\nabla q), \quad q\in\partial \mathcal E(u),\quad\text{in }\Omega\times (0,T),\\
& u(0,x)=u_0(x)\geq 0\quad\text{in }\Omega,
\end{aligned}
\end{equation}
with normalised mass $\int_\Omega u_0\; dx = 1$ and where the total variation (TV) functional $\mathcal E$ (see \cite{AFP,EvGa}) is defined by
\begin{equation}\label{totalvariation}
\mathcal E(u) := |Du|(\Omega)= \sup_{{\bf g} \in C_0^\infty(\Omega;\mathbb{R}^d), \norm{g}_\infty \leq 1} \int_\Omega  u ~\nabla \cdot {\bf g} ~dx,
\end{equation}
$d=1,2$. Equation \eqref{TVWassflow} can be formally derived as the $L^2$-Wasserstein gradient flow of the TV functional $\mathcal E$ in \eqref{totalvariation} and constitutes a nonlinear fourth-order diffusion equation. In this paper we study this equation as a regularising procedure for $u_0$ being a noisy image. Motivated by previous contributions in higher-order regularisation (see \cite{bredies,benning1,didas,kostas,setzer}), this approach promises to maintain the desirable properties of TV regularisation, such as preservation of edges in the image, while at the same time reducing well-known artifacts of TV such as \emph{staircaising}.

Originally, \eqref{TVWassflow} has been proposed in \cite{burg} for density estimation and smoothing. Therein, the authors propose to compute a smoothed version $u$ of a given probability density $u_0$ as a minimiser of
\begin{equation}\label{minprob}
\frac{1}{2} W_2(u_0{\mathcal L}^d,u{\mathcal L}^d)^2 + \alpha \mathcal E(u).
\end{equation}
 Here, $W_2(u_0{\mathcal L}^d,u {\mathcal L}^d)$ is the $L^2$--Wasserstein distance between $u_0{\mathcal L}^d$ and $u{\mathcal L}^d$ (${\mathcal L}^d$ denotes the Lebesgue measure in $\mathbb{R}^d$, $d=1,2$) and defines a distance within the space of probability measures \cite{gradflowbuch,Vi,Vi09,Ambrosioetal,Kant}. This minimisation problem can be interpreted as a discrete approximation of a solution of the gradient flow \eqref{TVWassflow} of $\mathcal E(u)$ with respect to the $L^2$-Wasserstein metric. More precisely, the minimisation of \eqref{minprob} represents one timestep of De Giorgi's minimising movement scheme (see, e.g. \cite{gradflowbuch,JKO}) to the functional $\mathcal E(u)$ with timestep $\alpha$. By construction the regularisation method \eqref{minprob} proposed in \cite{burg} is non-smooth, i.e., edge preserving, and conserves mass, i.e., is density preserving. In \cite{burg} the numerical solution of \eqref{minprob} has been accomplished by a combination of the Benamou-Brenier formulation \cite{benamoubrenier} for the Wasserstein distance, an augmented Lagrangian method, and an operator splitting technique \cite{Brune10}. This numerical procedure is in the flavour of several recently proposed numerical schemes for equations with gradient flow structure, cf., e.g., \cite{CM09,BCC08,DMM10,BCW10} and references therein.\\

Equation \eqref{TVWassflow} has been further investigated in \cite{dur}, where the authors numerically study the scale space properties and high-contrasting effects of the equation by solving it with a dimensional \emph{alternating direction implicit} (ADI) operator splitting approach. There, the subgradient $q$ of the TV seminorm in \eqref{TVWassflow} is approximated by an $\epsilon$-regularisation of the form
\begin{equation} \label{PDEwass}
q\approx \nabla\cdot \left(\frac{\nabla u}{\sqrt{|\nabla u|^2+\epsilon}} \right), \quad 1\gg \epsilon>0.
\end{equation}

From a computational point of view, finding numerical schemes that solve higher-order equations like \eqref{TVWassflow} is a challenging problem. Dealing with an evolutionary nonlinear fourth-order partial differential equation, we aim to find an efficient and reliable method avoiding a naive explicit discretisation in time that might present time step size restrictions (compare \cite{smereka}) and, because of the strong nonlinearity of the subgradients of TV, additionally add constraints to the stability condition of the discrete time stepping scheme, compare \cite{ChMu,BuDiWei,CDS}. \\

In this paper we propose a formulation of \eqref{TVWassflow} which characterises the elements of the subdifferential of TV in an alternative way. Instead of considering a characterisation of the type \eqref{PDEwass} for these elements, we use the approach proposed in \cite{benning,muller} and deal with a relaxed primal-dual formulation of \eqref{TVWassflow}. In \cite{CDS} the authors consider such a kind of approach applied to a similar fourth-order PDE as well as a directional splitting strategy that has proposed there as a direction for future research.

\paragraph{Notation.}
 We denote by $u$ the solution of the continuous equations and by $U$ the solution of the time discrete numerical schemes we are going to present. We write $U_n$ to indicate the approximation of $u(n\Delta t,\cdot),$ $n\geq 1$, where $\Delta t$ is the time step size. We will typically consider a rectangular domain $\Omega=[a,b]\times[c,d]$ and approximate it by a finite grid $\{a=x_1<\ldots<x_N=b\}\times\{c=y_1<\ldots<y_M=d\}$ with equidistant step-size $h=(b-a)/N=(d-c)/M$. For vectors, we will indicate their components using the superscripts notation: $\textbf{y}=\left( y^1\, y^2\right)^{\top}$. With a little abuse of notation, we will still use the notations $\nabla$ and $\nabla\cdot$ to indicate the differential operators discretised with either forward or backward finite differences applied to the discretised quantities (see Section \ref{sec:numres}).

\section{Primal-dual formulation of the $TV$-Wasserstein flow} \label{sec:primdual}
\setcounter{equation}{0}

We aim to characterise the elements in the subdifferential of the TV seminorm \eqref{totalvariation} by primal-dual iterations, as suggested, for instance, in \cite{benning,muller} when dealing with the classical second-order $TV$-denoising model. In what follows we discuss such a strategy combined with an implicit time stepping method for solving equation \eqref{TVWassflow}. By definition of the subdifferential, the property $q\in\partial|Du|(\Omega)$ means:
\begin{equation} \label{defsubdiffTV}
q\in \partial |Du|(\Omega)\iff |Du|(\Omega)-\int_\Omega qu \; dx \leq|Dv|(\Omega)-\int_\Omega qv \; dx ,\quad\forall v\in L^2(\Omega).
\end{equation}
Equivalently, if $u\in BV(\Omega)\subset L^2(\Omega)$ achieves the minimum of the following variational problem
\begin{equation} \label{minprobl}
\min_{u\in BV(\Omega)}\left\{|Du|(\Omega)-\int_\Omega qu\; dx\right\},
\end{equation}
then, by definition of being minimum, \eqref{defsubdiffTV} is fulfilled and then $q\in\partial|Du|(\Omega)$. Inserting the definition of the total variation seminorm \eqref{totalvariation} into \eqref{minprobl} we receive
\begin{equation} \label{minprob1}
\min_{u\in BV(\Omega)}\left\{\sup_{\textbf{p}\in C_0^{\infty}(\Omega;\R^2),\ \norm{\textbf{p}}_\infty\leq 1}\int_\Omega  u\nabla\cdot\textbf{p} \; dx -\int_\Omega qu \; dx  \right\}
\end{equation}
which is typically known as the \emph{primal-dual} formulation of the problem \eqref{TVWassflow}. The constraint on \textbf{p} appearing in \eqref{minprob1} can be relaxed, for instance, by a penalty method. To this end we remove the constraint from the minimisation in \eqref{minprob1} and instead add a term to the functional that penalises it if $\norm{\textbf{p}}_\infty>1$. A typical example for such a penalty term $F$ is
\begin{equation*}
F(s)=\frac{1}{2}\norm{\max\{s,0\}}^2_2.
\end{equation*}
With these considerations we reformulate \eqref{minprob1} into the following minimisation problem
\begin{equation} \label{minprob2}
\min_{u\in BV(\Omega)}\sup_{\textbf{p}\in C_0^\infty(\Omega;\R^2)}\left\{\int_\Omega u\nabla\cdot\textbf{p} \; dx -\frac{1}{\varepsilon}F(|\textbf{p}|-1)-\int_\Omega q u \; dx \right\},
\end{equation}
where the parameter $1\gg \varepsilon>0$ is small and measures the weight of our penalisation. We can then find the optimality conditions for both $\textbf{p}$ and $u$ in \eqref{minprob2} which, merged with the original equation \eqref{TVWassflow}, allow us to consider the following, alternative formulation of the $TV$-Wasserstein model:
\begin{equation}  \label{primdualsyst}
\left\{\begin{aligned}
 u_t & =\nabla\cdot(u\nabla q), \\
 q & =\nabla\cdot\textbf{p}, \\
 0 & = -\nabla u-\frac{1}{\varepsilon}H(\textbf{p}).
 \end{aligned}\right.
\end{equation}
In the system above $H$ denotes the derivative of the penalty term $F(|\textbf{p}|-1)$, i.e.
$$
H(\textbf p) = \mathbbm{1}_{\{|\textbf{p}|\geq 1\}} \mathrm{sgn}(\textbf p) (|\textbf{p}|-1),
$$
which we linearise via its first-order Taylor approximation
\begin{equation*}
H(\textbf p)\approx H(\tilde{\textbf p}) + H'(\tilde{\textbf p})(\textbf p-\tilde{\textbf p}),
\end{equation*}
where with $H'$ we indicate the Jacobian of $H$. In order to guarantee the invertibility of the now linear operator that defines the system, we add an additional damping term in $\mathbf p$, as suggested, for instance, in \cite{muller}. Collecting everything, we propose the following numerical scheme for solving \eqref{primdualsyst},
\begin{equation}  \label{primdualsystlin}
\left\{\begin{aligned}
  \frac{U^{(k)}_{n+1}-U_{n}}{\Delta t}&=\nabla\cdot(U_n\nabla Q^{(k)}_{n+1}) \\
  Q^{(k)}_{n+1}&=\nabla\cdot\textbf{P}^{(k)}_{n+1},\\
  0 &= -\nabla U^{(k)}_{n+1}-\frac{1}{\varepsilon}H(\textbf{P}^{(k-1)}_{n+1})\\
 &\quad -\frac{1}{\varepsilon}H'(\textbf{P}^{(k-1)}_{n+1})(\textbf{P}^{(k)}_{n+1}-\textbf{P}^{(k-1)}_{n+1})-\tau^k(\textbf{P}^{(k)}_{n+1}-\textbf{P}^{(k-1)}_{n+1}).
 \end{aligned}\right.
\end{equation}
We apply Newton's method to solve \eqref{primdualsystlin}. The scheme consists of two nested iterations. The subscripts $n$ are related to the outer time step evolution of the process evolving $U$. At each time step an implicit approximation of the quantities $U_{n+1}, Q_{n+1}$ and $\textbf{P}_{n+1}$ is obtained by the application of an inner damped Newton process that runs depending on the superscript $k$. The sequence of parameters $\tau^k$ controls the damping of the Newton iterations: it starts from a large value $\tau^0$ and then decreases, thus ensuring faster convergence. System \eqref{primdualsystlin} could now be discretised in space as described in Section \ref{sec:numres}. For computational simplicity we consider a slightly different penalty term $F$ (see \cite{muller}):
 \begin{equation*}
 F(\textbf{p})=F(p^1,p^2)=\frac{1}{2}\norm{\max\{|p^1|-1,0\}}^2_2+\frac{1}{2}\norm{\max\{|p^2|-1,0\}}^2_2,
 \end{equation*}
 that results into an anisotropic TV term. Whence
 \begin{equation*}
 H(p^1,p^2)=\left(\begin{aligned}
 & \mathrm{sgn}(p^1)(|p^1|-1)\mathbbm{1}_{\{|p^1|\geq 1\}}\\
 &\mathrm{sgn}(p^2)(|p^2|-1)\mathbbm{1}_{\{|p^2|\geq 1\}}
 \end{aligned}\right),\
 H'(p^1,p^2)=\left(\begin{aligned}
 & \mathbbm{1}_{\{|p^1|\geq 1\}} & 0 \\
 & 0 & \mathbbm{1}_{\{|p^2|\geq 1\}}
 \end{aligned}\right).
 \end{equation*}

\section{Numerical results}  \label{sec:numres}
\setcounter{equation}{0}

We discretise the gradient and the divergence operators appearing in the system above using standard forward and backward finite differences, thus preserving the adjointness properties of the operators as described in \cite{chamb} and used also in \cite{benning,muller}. We then build up the matrices representing system \eqref{primdualsystlin}. In each step of Newton's method the block-structure of the Jacobian matrix is exploited by inverting it with a Schur complement strategy. For all the following tests we use the following stopping criterion:
${\norm{U^{(k)}_{n+1}-U^{(k-1)}_{n+1}}}/{\norm{U^{(k)}_{n+1}}}\leq\epsilon_{tol},
$
where $\norm{\cdot}$ is the $\ell^2$-norm and $\epsilon_{tol}$ the tolerance.

\begin{figure}[!h]
\begin{center}
\includegraphics[height=3.5cm]{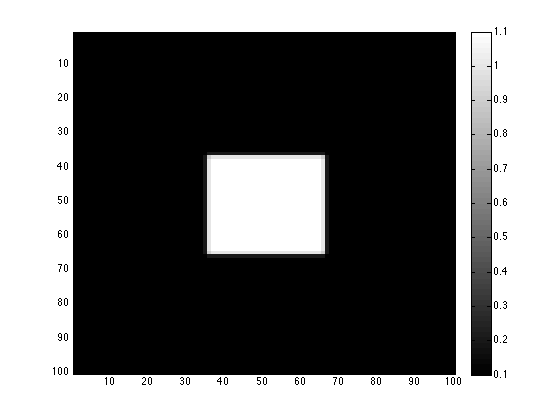}
\includegraphics[height=3.5cm]{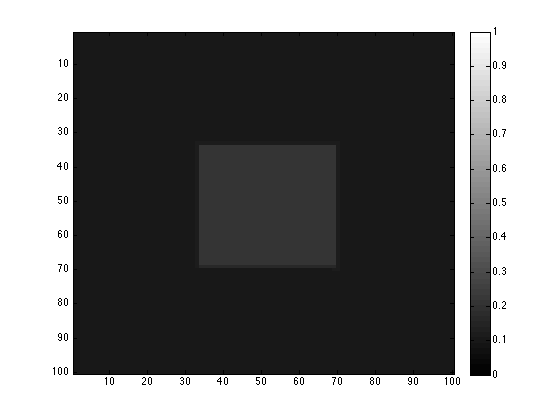}
\end{center}
\caption{Initial condition (l.) and solution of the $TV$-Wasserstein gradient flow after $10000$ iterations (r.). $\varepsilon=10^{-3}, \tau^0=1$.}
\label{square}
\end{figure}

Figure \ref{square} shows the solution of \eqref{TVWassflow} computed via \eqref{primdualsystlin} on a $100\times 100$ pixels initial condition of a square. The result shows the difference between the $TV$-Wasserstein approach and the pure $TV$ one. The latter, in fact, would just decrease the intensity of the square  and increase the intensity of the background without changing its support. Because of the mass preservation properties the solution of \eqref{TVWassflow} enlarge their support, as pointed out in \cite{burg} where the authors show the self-similarity of the solutions.

\begin{figure}[!h]
\begin{subfigure}{0.5\textwidth}
\centering
\includegraphics[height=3.8cm,width=4cm]{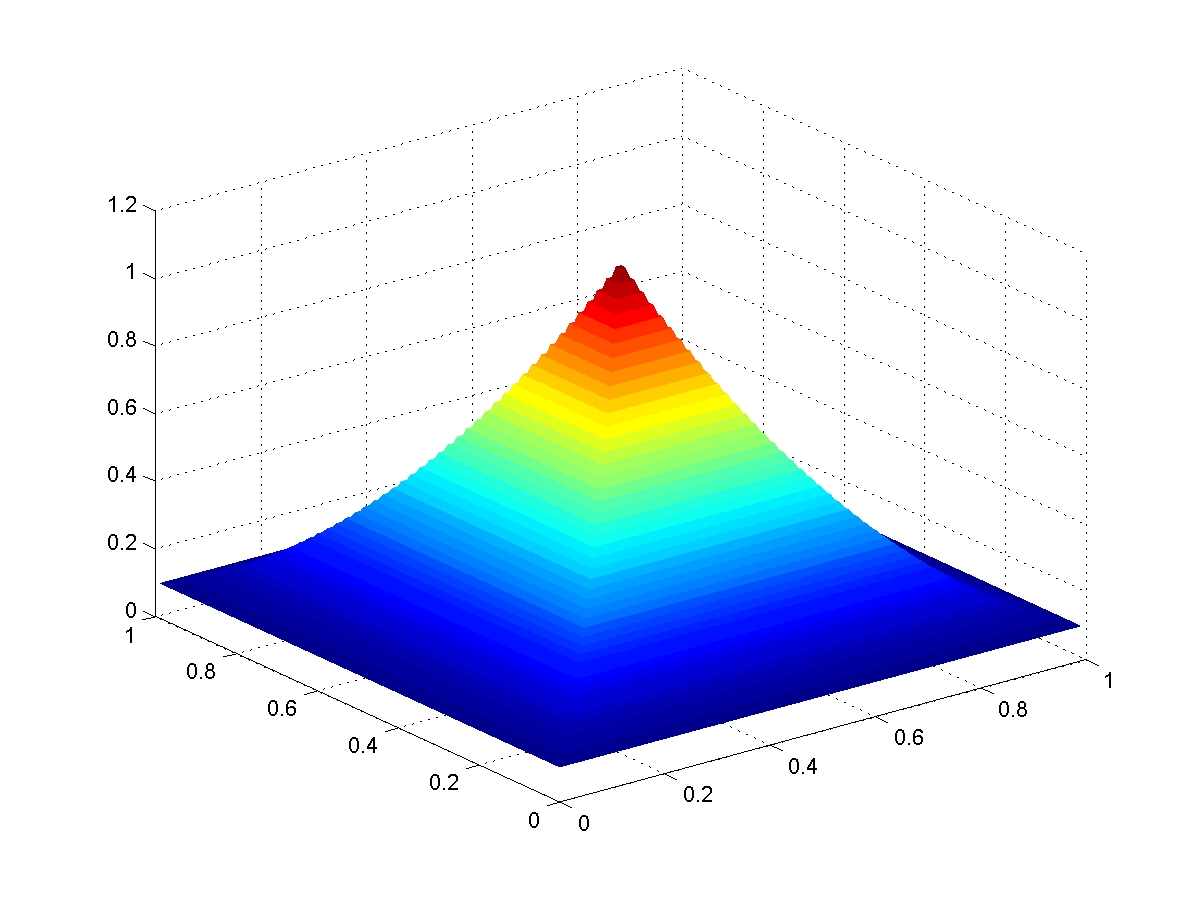}
\caption{Initial condition.}
\end{subfigure}
\quad
\begin{subfigure}{0.5\textwidth}
\centering
\includegraphics[height=3.8cm,width=4cm]{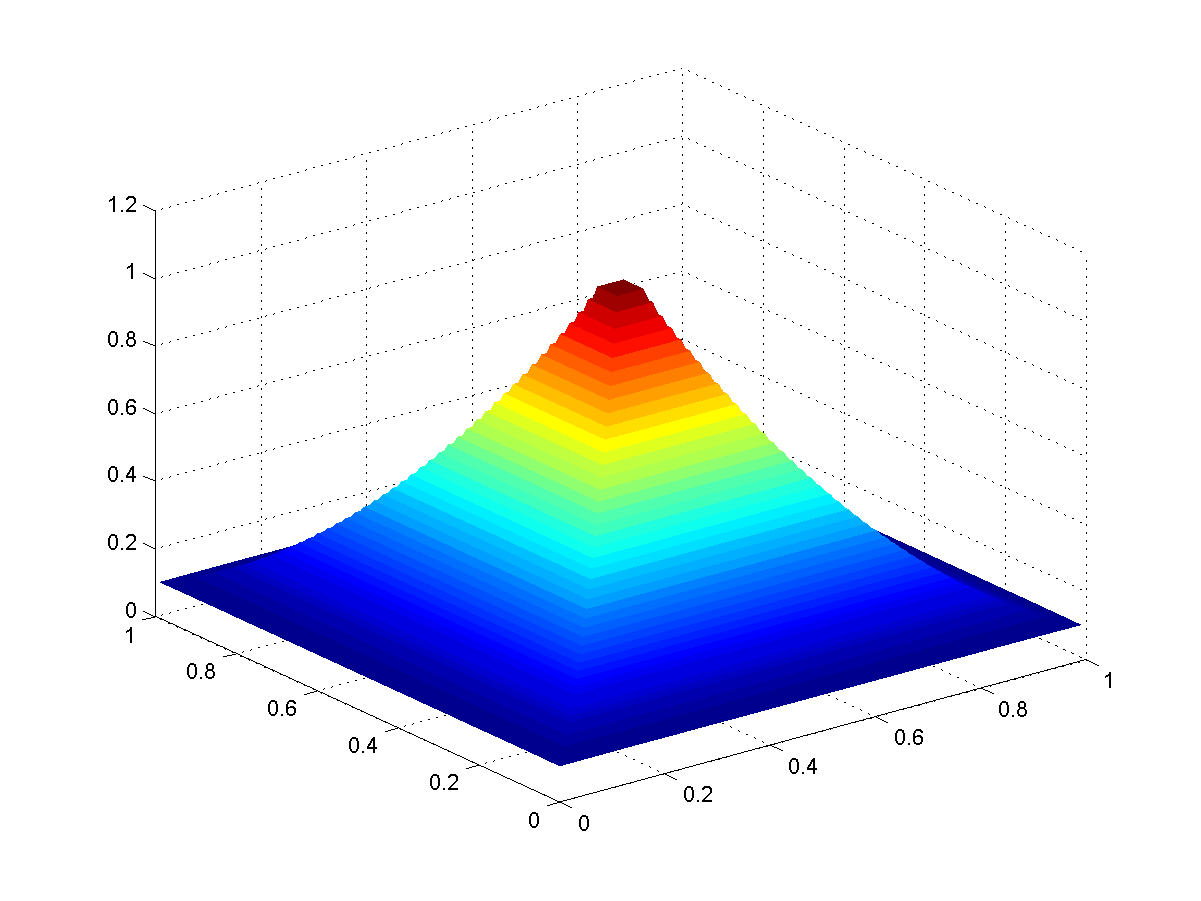}
\caption{Solution after $100$ iterations.}
\end{subfigure}\\
\begin{subfigure}{0.5\textwidth}
\centering
\includegraphics[height=3.8cm,width=4cm]{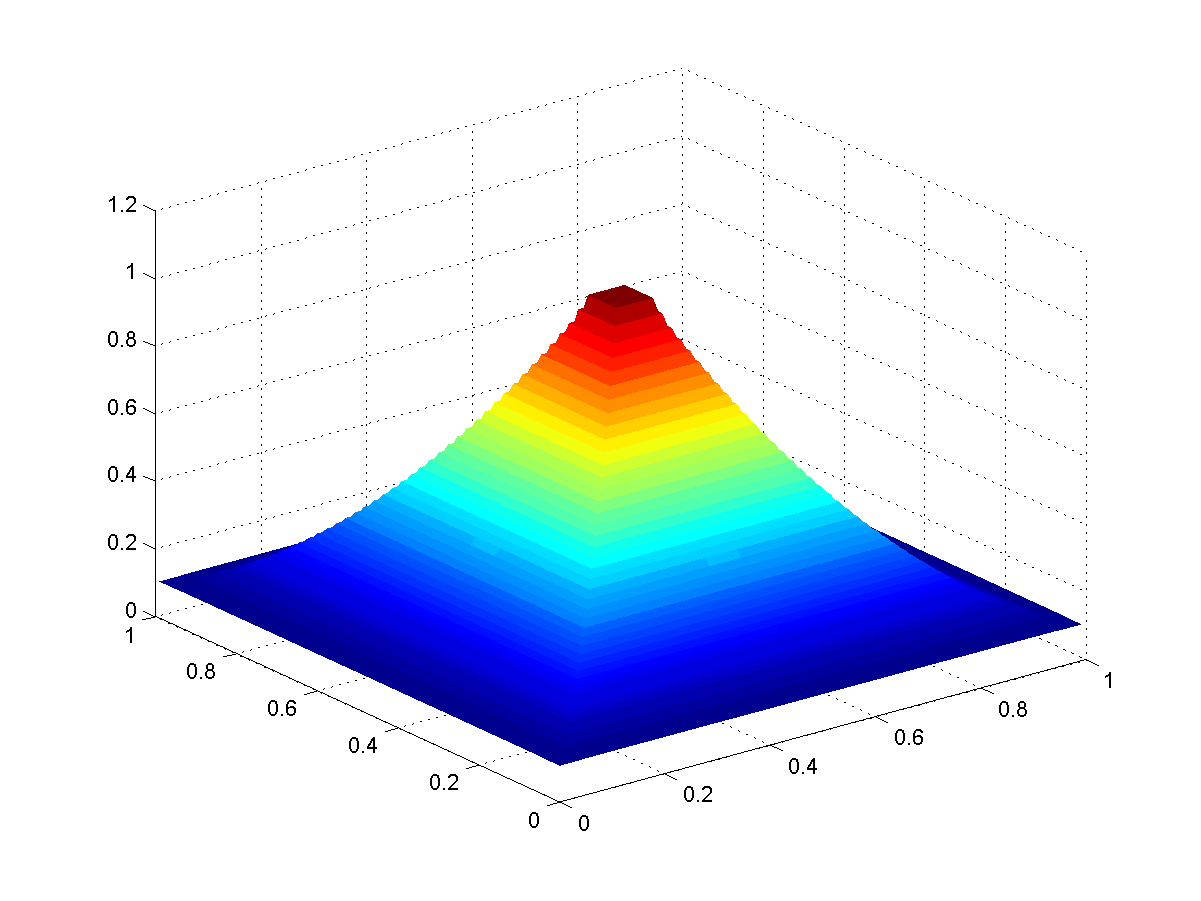}
\caption{Solution after $300$ iterations. }
\end{subfigure}
\quad
\begin{subfigure}{0.5\textwidth}
\centering
\includegraphics[height=3.8cm,width=4cm]{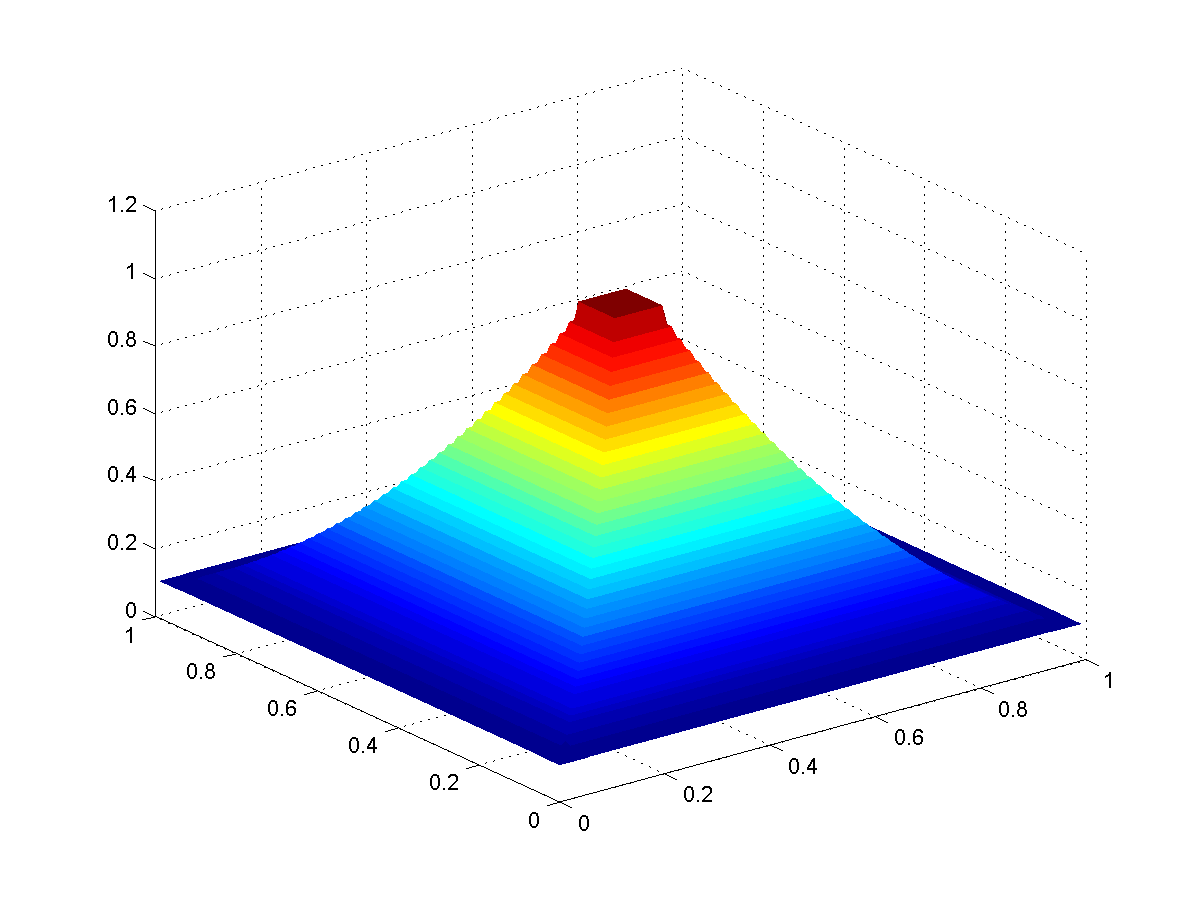}
\caption{Solution after $600$ iterations.}
\end{subfigure}
\caption{Time evolution of \eqref{primdualsystlin} with $\varepsilon=10^{-7}, \tau^0=1$.}
\label{pyramidtime}
\end{figure}

 Figure \ref{pyramidtime} shows the time evolution of the method \eqref{primdualsystlin} applied to a pyramidal initial condition defined on the domain $\Omega=[0,1]\times[0,1]$.


\begin{figure}[!h]
\begin{subfigure}{0.5\textwidth}
\centering
\includegraphics[height=3.8cm,width=4cm]{pyramidOriginal.png}
\caption{Initial condition.}
\end{subfigure}
\quad
\begin{subfigure}{0.5\textwidth}
\centering
\includegraphics[height=3.8cm,width=4cm]{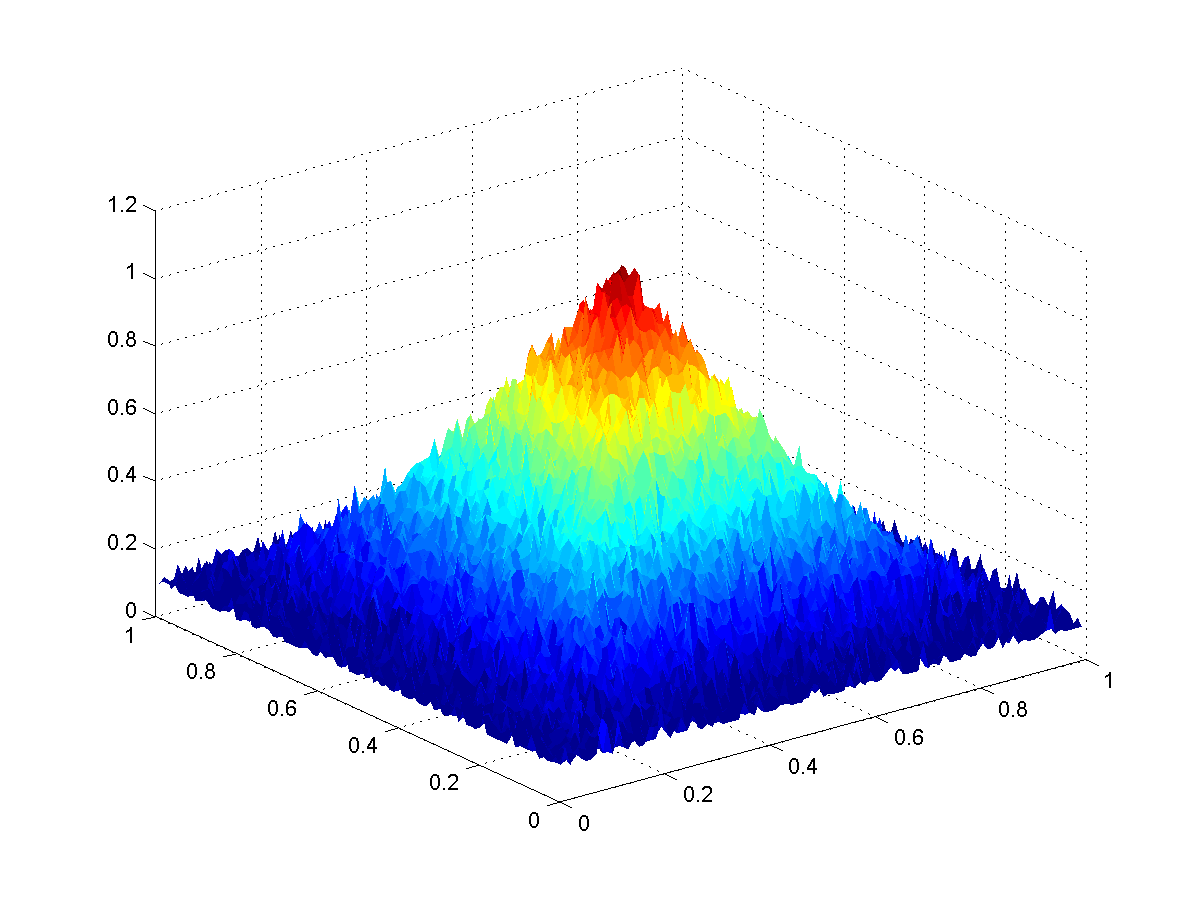}
\caption{Noisy pyramid with Gaussian noise of variance $0.001$.}
\end{subfigure}\\
\begin{subfigure}{0.5\textwidth}
\centering
\includegraphics[height=3.8cm,width=4cm]{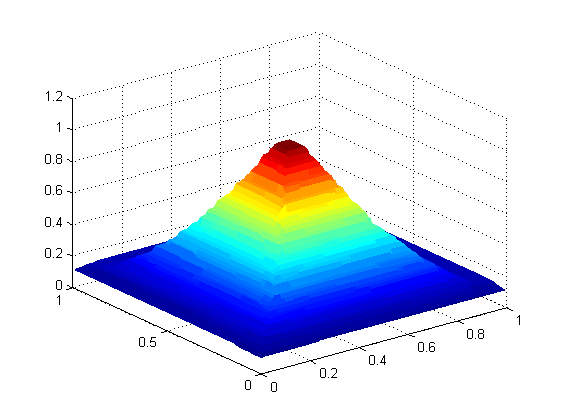}
\caption{Denoised pyramid with $TV$ primal-dual model. }
\end{subfigure}
\quad
\begin{subfigure}{0.5\textwidth}
\centering
\includegraphics[height=3.8cm,width=4cm]{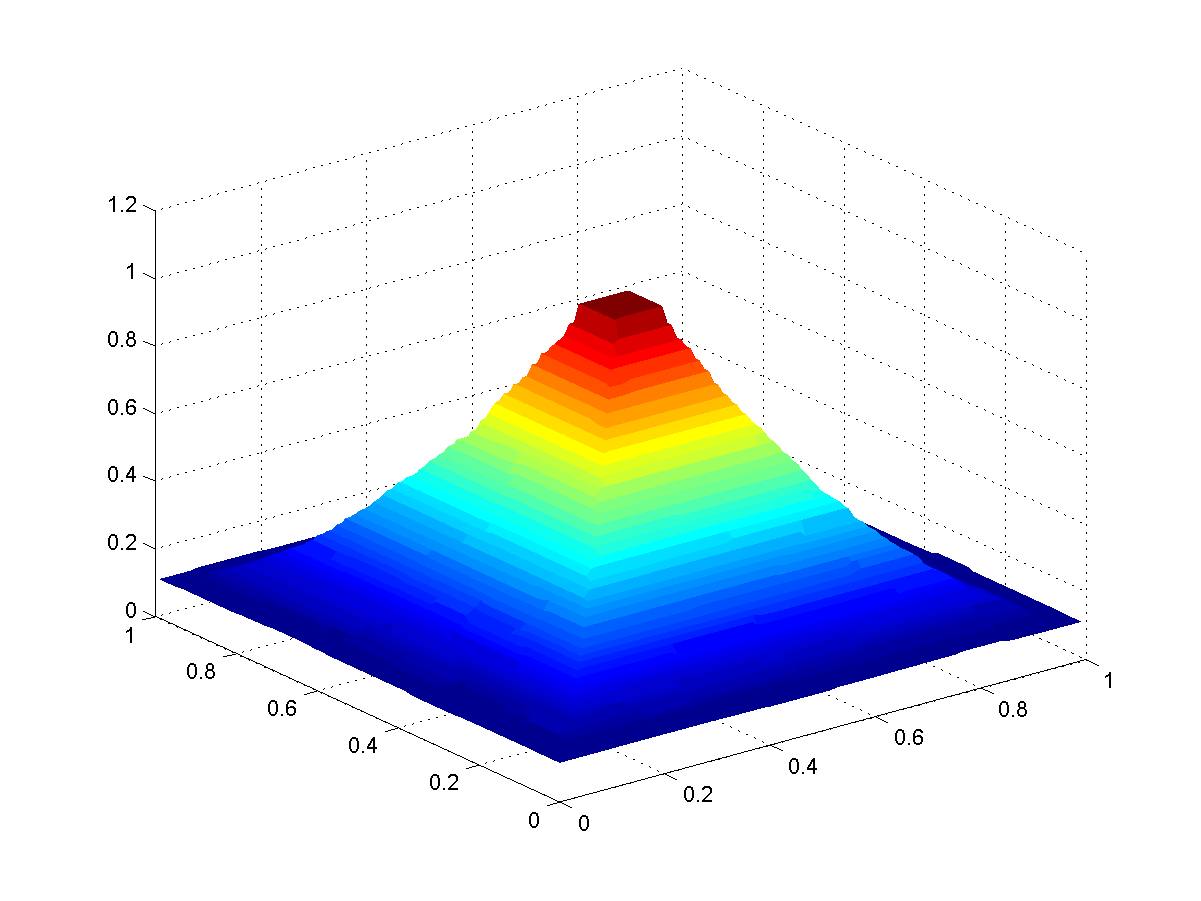}
\caption{Denoised pyramid with $TV$-Wasserstein primal-dual model \eqref{primdualsystlin}.}
\end{subfigure}
\caption{Denoising with $TV$ and $TV$-Wasserstein primal-dual methods. $\varepsilon=10^{-5}, \tau^0=1$.}
\label{pyramid}
\end{figure}

We now apply our method to a denoising problem. Figure \ref{pyramid} shows the pyramid and a noisy version of it obtained by adding Gaussian noise with variance $0.001$. The denoised version is obtained both with the primal-dual $TV$ method with penalty term described in \cite{muller} and with our method. We observe that while the simple application of the $TV$ model creates staircaising, the use of higher-order models reduces artifacts and preserves structures.

\begin{figure}[!h]
\begin{subfigure}{0.3\textwidth}
\centering
\includegraphics[height=3cm,width=3.5cm]{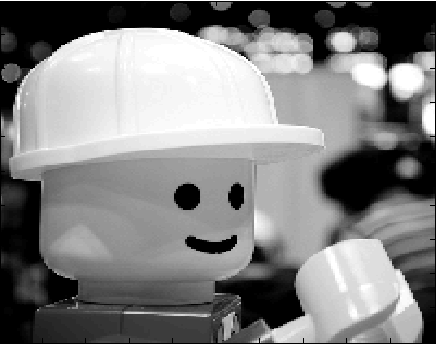}
\caption{Original.}
\end{subfigure}
\quad
\begin{subfigure}{0.3\textwidth}
\centering
\includegraphics[height=3cm,width=3.5cm]{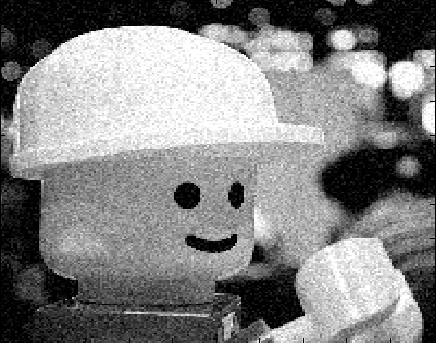}
\caption{Noisy image.}
\end{subfigure}
\quad
\begin{subfigure}{0.3\textwidth}
\centering
\includegraphics[height=3cm,width=3.5cm]{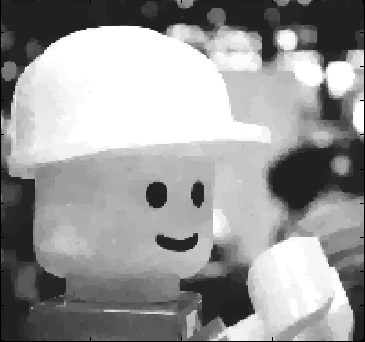}
\caption{$TV$ primal-dual result.}
\end{subfigure}\\
\begin{subfigure}{0.3\textwidth}
\centering
\includegraphics[height=3cm,width=3.5cm]{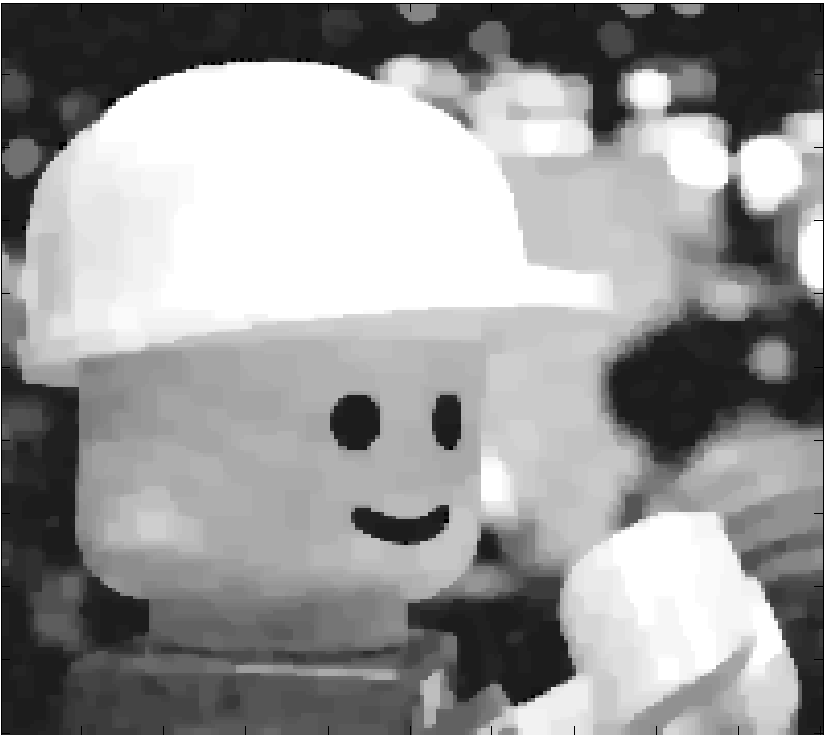}
\caption{$TV$-Wasserstein result after $100$ iterations.}
\end{subfigure}
\quad
\begin{subfigure}{0.3\textwidth}
\centering
\includegraphics[height=3cm,width=3.5cm]{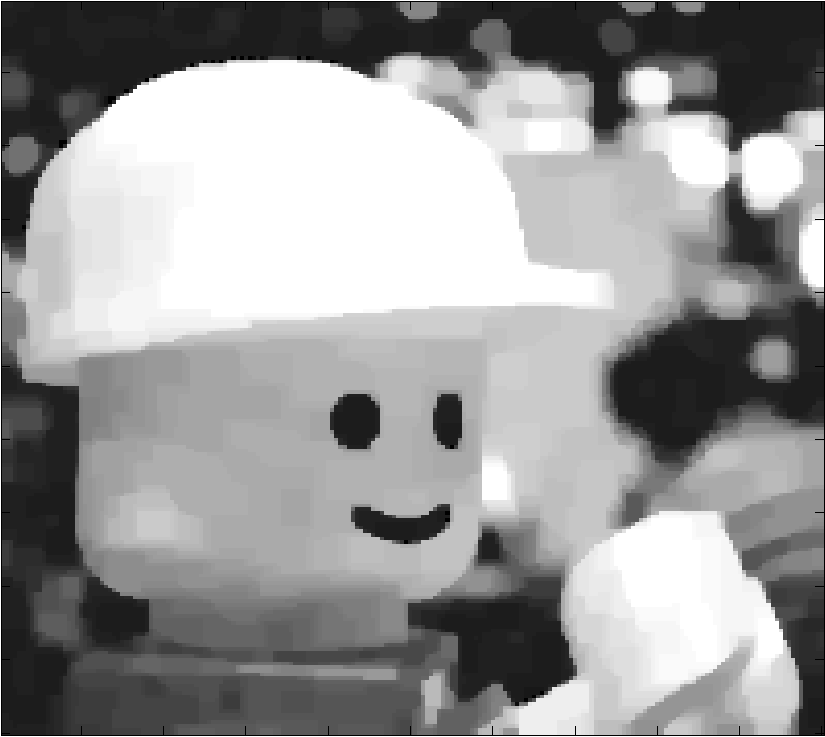}
\caption{$TV$-Wasserstein result after $200$ iterations.}
\end{subfigure}
\quad
\begin{subfigure}{0.3\textwidth}
\centering
\includegraphics[height=3cm,width=3.5cm]{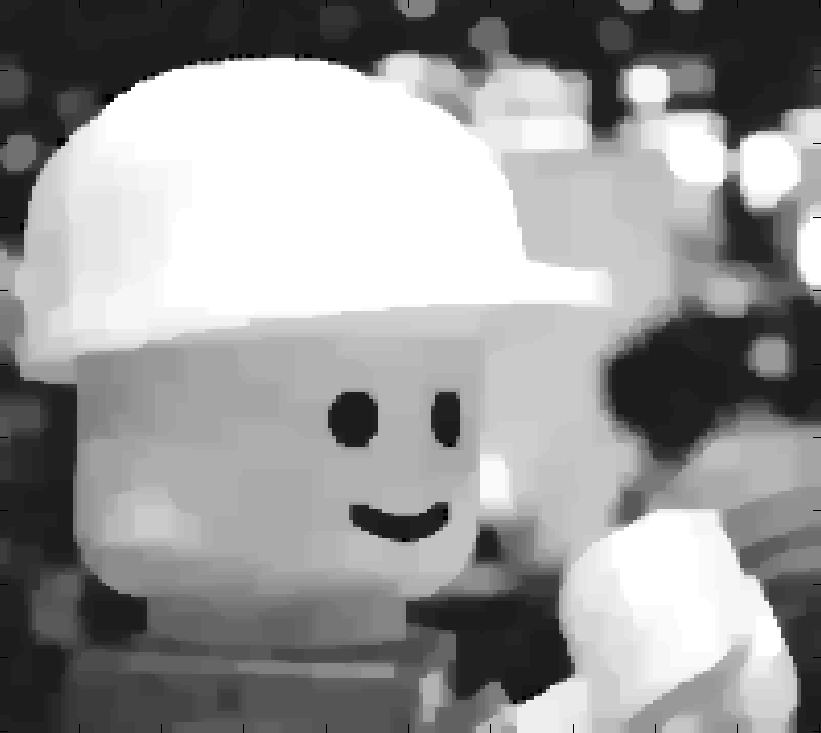}
\caption{$TV$-Wasserstein result after $350$ iterations.}
\end{subfigure}
\caption{Comparison with $TV$ primal-dual denoising model and $TV$-Wasserstein. $\varepsilon=10^{-7}, \tau^0=1$.}
\label{LEGO1}
\end{figure}

We consider in Figure \ref{LEGO1} a real-world image of a LEGO man. The dimension of the image is $200\times 200$ pixels. We add a Gaussian noise with zero mean and variance equal to $0.005$ and we show the result of time evolution of the process after some time iterations. A result with the application of the $TV$ primal-dual method is given for comparison as well.

\paragraph{Conclusions.} We proposed a numerical method to solve \eqref{TVWassflow}. Our strategy consists in a relaxation of the characterisation of the subgradients of the total variation via the addition of a \emph{penalty term}. The optimality conditions of the primal-dual formulation form a system of equations that can be solved via a damped Newton's method. The scheme shows good smoothing properties and reduced artifacts in comparison to the pure $TV$ method when applied to denoising problems.

\paragraph{Acknowledgements.}
Carola-Bibiane Sch\"{o}nlieb acknowledges financial support provided by the Cambridge Centre for Analysis (CCA), the Royal Society International Exchanges Award IE110314 for the project \emph{High-order Compressed Sensing for Medical Imaging}, the EPSRC first grant Nr. EP/J009539/1 \emph{Sparse \& Higher-order Image Restoration}, and the EPSRC / Isaac Newton Trust Small Grant on \emph{Non-smooth geometric reconstruction for high resolution MRI imaging of fluid transport in bed reactors}. Further, this publication is based on work supported by Award No.
KUK-I1-007-43, made by King Abdullah University of Science and Technology (KAUST).

\end{document}